\documentclass[a4paper,12pt]{amsart}
\usepackage{amssymb}
\usepackage{ifthen}
\usepackage[dvips]{graphicx}
\nonstopmode \numberwithin{equation}{section}
\usepackage{amssymb}
\usepackage{ifthen}
\usepackage{graphicx}
\usepackage{amsmath}
\usepackage[T1]{fontenc} %skandit

\nonstopmode \numberwithin{equation}{section}
\theoremstyle{plain}
\newtheorem{thm}[equation]{Theorem}
\newtheorem{cor}[equation]{Corollary}
\newtheorem{lem}[equation]{Lemma}
\newtheorem{prop}{Proposition}

\newtheorem{conj}{Conjecture}
\newtheorem{innercustomthm}{Theorem}
\newenvironment{customthm}[1]
{\renewcommand\theinnercustomthm{#1}\innercustomthm}
{\endinnercustomthm}

\theoremstyle{definition}
\newtheorem{defn}{Definition}[section]

\newtheorem{prob}{Problem}
\newtheorem{rem}{Remark}[section]

\newcounter{minutes}
\divide\time by 60
\newcounter{hours}
\multiply\time by 60
\addtocounter{minutes}{-\time}
\newcounter {own}
\def\theown {\thesection       .\arabic{own}}

\newenvironment{pf}[1][]{%
	\vskip 3mm
	\noindent
	\ifthenelse{\equal{#1}{}}%
	{{\slshape Proof. }}%
	{{\slshape #1.} }%
}%
{\qed\bigskip}

\newcounter{alphabet}
\newcounter{tmp}

\newcommand{\real}{{\operatorname{Re}\,}}

\def\be{\begin{equation}}
\def\ee{\end{equation}}

\newcommand{\bee}{\begin{enumerate}}
	\newcommand{\eee}{\end{enumerate}}

\newcommand{\blem}{\begin{lem}}
	\newcommand{\elem}{\end{lem}}
\newcommand{\bthm}{\begin{thm}}
	\newcommand{\ethm}{\end{thm}}
\newcommand{\bcor}{\begin{cor}}
	\newcommand{\ecor}{\end{cor}}
\newcommand{\beg}{\begin{examp}}
	\newcommand{\eeg}{\end{examp}}
\newcommand{\begs}{\begin{examples}}
	\newcommand{\eegs}{\end{examples}}

\newcommand{\bdefn}{\begin{defn}}
	\newcommand{\edefn}{\end{defn}}

\newcommand{\bprob}{\begin{prob}}
	\newcommand{\eprob}{\end{prob}}
\newcommand{\bei}{\begin{itemize}}
	\newcommand{\eei}{\end{itemize}}

\newcommand{\bcon}{\begin{conj}}
	\newcommand{\econ}{\end{conj}}
\newcommand{\bcons}{\begin{conjs}}
	\newcommand{\econs}{\end{conjs}}
\newcommand{\bprop}{\begin{prop}}
	\newcommand{\eprop}{\end{prop}}
\newcommand{\br}{\begin{rem}}
	\newcommand{\er}{\end{rem}}
\newcommand{\brs}{\begin{rems}}
	\newcommand{\ers}{\end{rems}}
\newcommand{\bo}{\begin{obser}}
	\newcommand{\eo}{\end{obser}}
\newcommand{\bos}{\begin{obsers}}
	\newcommand{\eos}{\end{obsers}}
\newcommand{\bpf}{\begin{pf}}
	\newcommand{\epf}{\end{pf}}
\newcommand{\ba}{\begin{array}}
	\newcommand{\ea}{\end{array}}
\newcommand{\beq}{\begin{eqnarray}}
\newcommand{\beqq}{\begin{eqnarray*}}
\newcommand{\eeq}{\end{eqnarray}}
\newcommand{\eeqq}{\end{eqnarray*}}

\begin{document}

\title{Bohr radius for certain classes of starlike and convex univalent functions}

\author{Vasudevarao Allu}
\address{Vasudevarao Allu,
	School of Basic Science,
	Indian Institute of Technology Bhubaneswar,
	Bhubaneswar-752050, Odisha, India.}
\email{avrao@iitbbs.ac.in}

\author{Himadri Halder}
\address{Himadri Halder,
	School of Basic Science,
	Indian Institute of Technology Bhubaneswar,
	Bhubaneswar-752050, Odisha, India.}
\email{hh11@iitbbs.ac.in}

\subjclass[{AMS} Subject Classification:]{Primary 30C45, 30C50, 30C80}

\def\thefootnote{}
\footnotetext{ {\tiny File:~\jobname.tex,
		printed: \number\year-\number\month-\number\day,
		\thehours.\ifnum\theminutes<10{0}\fi\theminutes }
} \makeatletter\def\thefootnote{\@arabic\c@footnote}\makeatother

\begin{abstract}
We say that a class $\mathcal{F}$ consisting of analytic functions $f(z)=\sum_{n=0}^{\infty} a_{n}z^{n}$ in the unit disk $\mathbb{D}:=\{z\in \mathbb{C}: |z|<1\}$ 
satisfies a Bohr phenomenon if there exists $r_{f} \in (0,1)$ such that 
$$
 \sum_{n=1}^{\infty} |a_{n}z^{n}|\leq d(f(0),\partial f(\mathbb{D}))
$$
for every function $f \in \mathcal{F}$ and $|z|=r\leq r_{f}$, where $d$ is the Euclidean distance. The largest radius $r_{f}$ is the Bohr radius for the class $\mathcal{F}$. 
In this paper, we establish the Bohr phenomenon for the classes consisting of Ma-Minda type starlike functions and Ma-Minda type convex functions as well as for the class of starlike functions with respect to a boundary point. 
%	The Bohr phenomenon for analytic functions of the form $f(z)=\sum_{n=0}^{\infty} a_{n}z^{n}$, first introduced by Harald Bohr in 1914, deals with finding the largest radius 
%	$r_{f}$, $0<r_{f}<1$, such that the inequality $\sum_{n=0}^{\infty} |a_{n}z^{n}| \leq 1$ holds whenever the inequality $|f(z)|\leq 1 $ holds in the unit disk 
%	$\mathbb{D}=\{z \in \mathbb{C}: |z|<1\}$. The exact value of this largest radius known as Bohr radius, which has been established to be  $r_{f}=1/3$. 
%	The Bohr phenomenon \cite{Abu-2010} for harmonic functions  $f$ of the form  $f(z)=h(z)+\overline {g(z)}$, where $h(z)=\sum_{n=0}^{\infty} a_{n}z^{n}$ 
%	and $g(z)=\sum_{n=1}^{\infty} b_{n}z^{n}$ is to find the largest radius $r_{f}$, $0<r_{f}<1$ such that 
%	$$\sum\limits_{n=1}^{\infty} (|a_{n}|+|b_{n}|) |z|^{n}\leq d(f(0),\partial f(\mathbb{D})) 
%	%\quad\mbox { for } |z|\leq r_{f}.
%	$$ 
%	holds for $|z|\leq r_{f}$, here $d(f(0),\partial f(\mathbb{D})) $ denotes the Euclidean distance between $f(0)$ and the boundary of $f(\mathbb{D})$. 
%	In this paper, we investigate the Bohr radius for several classes of harmonic functions in the unit disk $\mathbb{D}.$
\end{abstract}

\maketitle
\pagestyle{myheadings}
\markboth{Vasudevarao Allu and  Himadri Halder}{Bohr radius for certain classes of starlike and convex univalent functions}

\section{Introduction and Preliminaries}

Let $\mathcal{H}$ denote the class of all analytic functions $f$ in the unit disk $\mathbb{D}:=\{z\in \mathbb{C}: |z|<1\}$. Each function $f\in \mathcal{H}$ has the following power series  representation
\begin{equation} \label{him-p2-e-1.1}
f(z)=\sum_{n=0}^{\infty} a_{n}z^{n}.
\end{equation}
For $f\in \mathcal{H}$, the majorant series is defined by $M_{f}(r):= \sum_{n=0}^{\infty} |a_{n}|r^{n}$ for $|z|=r<1$.
 In 1914, Harald Bohr \cite{Bohr-1914} obtained the following remarkable result.
\begin{customthm}{A}
Let $f\in \mathcal{H}$ be given by \eqref{him-p2-e-1.1} and $|f(z)|<1$ for all $z\in \mathbb{D}$. Then 
\begin{equation} \label{him-p2-e-1.2}
M_{f}(r):= \sum_{n=0}^{\infty} |a_{n}|r^{n} \leq 1
\end{equation}
for all $z \in \mathbb{D}$ with $|z|=r \leq 1/3$. The constant $1/3$, called the Bohr radius, cannot be improved.
\end{customthm}
Harald Bohr initially obtained the inequality \eqref{him-p2-e-1.2} for $r\leq 1/6$. Later, Weiner, Riesz and Schur have 
independently established the inequality \eqref{him-p2-e-1.2} for $r\leq 1/3$ and have shown that the constant $1/3$ cannot be improved (see \cite{paulsen-2002,sidon-1927,tomic-1962}). The inequality \eqref{him-p2-e-1.2} is popularly known as Bohr's inequality. The Bohr's inequality has been extended to several complex variables by finding the multidimensional Bohr radius (see \cite{aizn-2000,aizenberg-2001,aizn-2007,boas-1997,boas-2000}). 

\vspace{4mm}
The Bohr inequality \eqref{him-p2-e-1.2} can be written in the following equivalent form 
\begin{equation} \label{him-p2-e-1.3}
d\left(\sum_{n=0}^{\infty} |a_{n}z^{n}|, |a_{0}|\right)= \sum_{n=1}^{\infty} |a_{n}z^{n}|\leq d(f(0),\partial f(\mathbb{D}))
\end{equation}
for $|z|=r \leq 1/3$, where $d$ is the Euclidean distance. We say that $\mathcal{H}$ satisfies the Bohr phenomenon if every function $f \in \mathcal{H}$ satisfies the inequality \eqref{him-p2-e-1.3}.
%If every function $f \in \mathcal{H}$ satisfies the inequality \eqref{him-p2-e-1.3}, then sometimes $\mathcal{H}$ is said to satisfy the classical Bohr phenomenon. 
The above equivalent form \eqref{him-p2-e-1.3} makes  the notion of the Bohr phenomenon evident for a class $\mathcal{F}$ consisting of analytic functions of the form \eqref{him-p2-e-1.1} which map the unit disk $\mathbb{D}$ into a given domain $\varOmega \subseteq \mathbb{C}$ such that $f(\mathbb{D}) \subseteq \varOmega$. The class $\mathcal{F}$ is said to satisfy the Bohr phenomenon if there exists $r_{\varOmega} \in (0,1)$ such that every function $f \in \mathcal{F}$ satisfies the inequality 
\begin{equation} \label{him-p2-e-1.4}
d\left(\sum_{n=0}^{\infty} |a_{n}z^{n}|, |a_{0}|\right)= \sum_{n=1}^{\infty} |a_{n}z^{n}|\leq d(f(0),\partial f(\mathbb{D}))
\end{equation}
for all $|z|=r \leq r_{\varOmega}$. The largest radius $r_{\varOmega}$ is called the Bohr radius for the class $\mathcal{F}$. The question arises as to whether or not the largest radius $r_{\varOmega}$ is changed with respect to different types of domains. The answer is affirmative. For any proper simply connected domain $\varOmega$, Abu-Muhanna \cite{Abu-2010} has proved that the sharp radius is $r_{\varOmega}=3-2\sqrt{2}$ for the class $\mathcal{F}$. For a convex domain $\varOmega$, Aizenberg \cite{aizn-2007} has shown that $r_{\varOmega}$ coincides with the classical Bohr radius 1/3. Bohr phenomenon for harmonic mappings has also been extensively studied by several authors (see \cite{abu-2014,Himadri-Vasu-P1,evdoridis-2019,kayumov-2018-a,kayumov-2018-b,kayumova-2018,Liu-2019}). For more information about Bohr phenomenon stated above and further related intriguing aspects, we refer the reader to the articles (see \cite{abu-2011,abu-2013,Ali-2017,Ali-2019,bene-2004,kayumov-2019}).

\vspace{5mm}
Let $\mathcal{A}$ denote the subclass of $\mathcal{H}$ consisting of functions normalized by $f(0)=f'(0)-1=0$. Each $f \in \mathcal{A}$ has the following form 
\begin{equation} \label{him-p2-e-1.5}
f(z)=z+ \sum_{n=2}^{\infty} a_{n}z^{n}.
\end{equation}
Let $\mathcal{S} \subseteq \mathcal{A}$ be the family of univalent ({\it i.e.} one-to-one) functions. Let $\mathcal{S}^{*}(\alpha)$ and $\mathcal{C}(\alpha)$ be the subclasses of $\mathcal{S}$ consisting 
of functions starlike of order $\alpha \, (0 \leq \alpha <1)$ and convex functions of order $\alpha \, (0 \leq \alpha <1)$ respectively. A function $f\in\mathcal{S}$ belongs to  $\mathcal{S}^{*}(\alpha)$  (respectively $\mathcal{C}(\alpha)$) if 
$\real({zf'(z)}/{f(z)})>\alpha$  for $z\in\mathbb{D}$   $(\real(1+{zf''(z)}/{f'(z)})>\alpha$ for $z\in\mathbb{D}$ {respectively}).
 It is known that $f \in \mathcal{C}(\alpha)$ if, and only if, $zf' \in \mathcal{S}^{*}(\alpha)$. The classes $\mathcal{S}^{*}:=\mathcal{S}^{*}(0)$ and $\mathcal{C}:=\mathcal{C}(0)$ are the family of starlike and convex functions in $\mathbb{D}$ respectively. For more properties of starlike and convex functions, we refer the reader to \cite{duren-1983,graham-2003,vasu-book}. For two analytic functions $f$ and $g$ in $\mathbb{D}$, we say that $f$ is subordinate to $g$, written $f \prec g$ in $\mathbb{D}$, if there exists a schwarzian function $\omega :\mathbb{D} \rightarrow \mathbb{D}$ with $\omega(0)=0$ such that $f(z)=g(\omega(z))$ for $z \in \mathbb{D}$. In particular, if $g$ is univalent in $\mathbb{D}$, then $f\prec g$ if, and only if, $f(0)=g(0)$ and $f(\mathbb{D}) \subseteq g(\mathbb{D})$. In 1992, Ma-Minda \cite{ma minda-1992-a} introduced the function classes $\mathcal{S}^{*}(\phi)$ and $\mathcal{C}(\phi)$ by unifying several earlier results on subordination. For the brevity, we write.
\begin{defn} \label{him-p2-def-1.1}
Let $\mathcal{S}^{*}(\phi)$ and $\mathcal{C}(\phi)$ denote the subclasses of $\mathcal{S}$ consisting of functions characterized by
$$\dfrac{zf'(z)}{f(z)} \prec \phi (z) \quad \mbox{and} \quad 1+ \dfrac{zf''(z)}{f'(z)} \prec \phi (z)$$
 respectively,
%$$
%\mathcal{S}^{*}(\phi):= \Biggl\{f \in \mathcal{S} : \dfrac{zf'(z)}{f(z)} \prec \phi (z)\Biggr\} \quad \mbox{and } \quad \mathcal{C}(\phi):= \Biggl\{f \in \mathcal{S} :1+ \dfrac{zf''(z)}{f'(z)} \prec \phi (z)\Biggr\},
%$$
where $\phi : \mathbb{D} \rightarrow \mathbb{D}$ is called Ma-Minda function which is analytic and univalent in $\mathbb{D}$ such that $\phi(\mathbb{D})$ has positive real part, symmetric with respect to the real axis, starlike with respect to $\phi(0)=1$ and $\phi ' (0)>0$.  Let the Taylor series expansion of $\phi(z)$ be of the form 
\begin{equation}  \label{him-p2-e-1.6}
\phi(z)=1+ \sum_{n=1}^{\infty} B_{n}z^{n} \quad (B_{1}>0)
\end{equation}
for $z \in \mathbb{D}$. Let $\mathcal{M}$ denote the class of all Ma-Minda functions in $\mathbb{D}$.
\end{defn}
We call $\mathcal{S}^{*}(\phi)$ and $\mathcal{C}(\phi)$ are the Ma-Minda type starlike and Ma-Minda type convex classes associated with $\phi$.
It is known that $f \in \mathcal{C}(\phi)$ if, and only if, $zf' \in \mathcal{S}^{*}(\phi)$. 
%Let the Taylor series expansion of $\phi(z)$ be of the form 
%\begin{equation}  \label{him-p2-e-1.6}
%\phi(z)=1+ \sum_{n=1}^{\infty} B_{n}z^{n} \quad (B_{1}>0)
%\end{equation}
%for $z \in \mathbb{D}$. 
For $\phi(z)=(1+z)/(1-z)$ we have $\mathcal{C}(\phi):=\mathcal{C}$ and $\mathcal{S}^{*}(\phi):=\mathcal{S}^{*}$. Evidently, for every such $\phi$ described in  Definition \ref{him-p2-def-1.1}, $\mathcal{S}^{*}(\phi)$ and $\mathcal{C}(\phi)$ are always subclasses of the classes $\mathcal{S}^{*}$ and $\mathcal{C}(\phi)$ respectively.

\vspace{4mm}
It is worth noting that for particular choices of $\phi$, the classes $\mathcal{S}^{*}(\phi)$ and $\mathcal{C}(\phi)$ generate several important subclasses of starlike and convex functions, respectively. For instance, $\mathcal{S}^{*}(\alpha):=\mathcal{S}^{*}\left((1+(1-2\alpha))/(1-z)\right)$ and 
$\mathcal{C}(\alpha):=\mathcal{C}\left((1+(1-2\alpha))/(1-z)\right)$, the Janowski starlike class $\mathcal{S}^{*}[A,B]:=\mathcal{S}^{*}((1+Az)/(1+Bz))$ and Janowski convex class $\mathcal{C}[A,B]:=\mathcal{C}((1+Az)/(1+Bz))$, where $-1\leq B < A \leq 1$. If $\phi(z)=\left((1+z)/(1-z)\right)^{\alpha}$ for $0< \alpha \leq 1$, then $\mathcal{C}(\alpha)$ and $\mathcal{S}^{*}(\alpha)$ are the classes of strongly convex and strongly starlike functions of order $\alpha$ (see \cite{ma minda-1992-b}). For
\begin{equation} \label{him-p2-e-1.7}
\phi(z)=1+ \dfrac{2}{\pi ^{2}}\left(log \dfrac{1+\sqrt{z}}{1-\sqrt{z}}\right)^{2},
\end{equation}
the class $\mathcal{C}(\phi) \, (\mathcal{S}^{*}(\phi) \quad \mbox{respectively})$ is the family {\it{ UCV}} ({\it UST}  respectively) of normalized uniformly convex (starlike respectively) functions introduced by Goodman (see \cite{goodman-1991-a,goodman-1991-b,ronning-1991,ronning-1993}). Ma and Minda (see \cite{ma minda-1992-b,ma minda-1993}) have studied extensively the class {\it UCV}. 
For $0 \leq \alpha <1$, Khatter {\it et al.} \cite{khatter-2019} introduced $\mathcal{S}^{*}_{\alpha,e}:=\mathcal{S}^{*}(\alpha +(1-\alpha)e^{z})$. For $\alpha=0$, $\mathcal{S}^{*}_{\alpha,e}$ 
reduces to $\mathcal{S}^{*}_{e}:=\mathcal{S}^{*}(e^{z})$ (see \cite{mendiratta-2015}). 
When $\phi(z)=1+4z/3+2z^{2}/3$, the unit disk $\mathbb{D}$ is mapped onto a domain bounded by a $cardiod$ and corresponding Ma-Minda starlike class $\mathcal{S}^{*}(\phi)$ reduces to the class $\mathcal{S}^{*}_{C}$ (see  \cite{sharma-2016}). For $\phi(z)=1+z/(1-\alpha z^{2})$, Kargar {\it et al.}  \cite{kargar-2019} have  introduced the class
$\mathcal{BS}^{*}(\alpha):= \mathcal{S}^{*}(1+z/(1-\alpha z^{2}))$, which is  associated with the {\it Booth lemniscate}. In 2016, Kumar and Ravichandran \cite{kumar-2016} considered the class $\mathcal{S}^{*}_{R}:=\mathcal{S}^{*}(\phi_{0})$, where $\phi_{0}$ is the rational function 
\begin{equation} \label{him-p2-e-1.8}
\phi_{0}(z)=1+\dfrac{z}{k} \left(\dfrac{k+z}{k-z}\right), \quad k=\sqrt{2} +1.
\end{equation}
For $\phi(z)=(1+sz)^{2}$ with $0< s \leq 1/\sqrt{2}$, the class $\mathcal{S}^{*}(\phi)$ reduces to $\mathcal{ST}_{L}(s):=\mathcal{S}^{*}\left(\left(1+sz\right)^{2}\right)$ (see  \cite{masih-2020}).

\vspace{5mm}
Ma and Minda have defined the following functions $h, k\in\mathcal{S}$ by
\begin{equation} \label{him-p2-e-1.9}
\dfrac{zh'(z)}{h(z)}=\phi (z) \quad \mbox{and} \quad 1+ \dfrac{zk''(z)}{k'(z)}=\phi (z).
\end{equation}
Here $h$ and $k$  play the role of Koebe function for the classes $\mathcal{S}^{*}(\phi)$ and $\mathcal{C}(\phi)$,  respectively.
Clearly, $h$ and $k$ belong to $\mathcal{S}^{*}(\phi)$ and $\mathcal{C}(\phi)$ respectively and $z k'(z)=h(z)$. The subordination and growth estimate for the class $\mathcal{S}^{*}(\phi)$ have been obtained by Ma and Minda (see \cite{ma minda-1992-a}).
\begin{lem} \label{him-p2-lem-1.10} \cite{ma minda-1992-a}
Let $f \in \mathcal{S}^{*}(\phi)$. Then $zf'(z)/f(z) \prec zh'(z)/h(z)$ and $f(z)/z \prec h(z)/z$.
\end{lem}
\begin{lem} \label{him-p2-lem-1.11} \cite{ma minda-1992-a}
Assume $f \in \mathcal{S}^{*}(\phi)$ and $|z|=r<1$. Then 
\begin{equation} \label{him-p2-e-1.12}
-h(-r) \leq |f(z)| \leq h(r).
\end{equation} 
Equality holds for some $z \neq 0$ if, and only if, f is a rotation of $h$.
\end{lem}
It has been pointed out in \cite{ma minda-1992-a} that $-h(-r)$ is increasing in $(0,1)$ and bounded by $1$  because  each $f \in \mathcal{S}^{*}(\phi)$ is normalized by $f(0)=f'(0)-1=0$.
Therefore,  $\lim_{r \rightarrow 1} -h(-r)$ exists and denote it by $-h(-1)$.
The following subordination and growth theorem for the class $\mathcal{C}(\phi)$ have been established in \cite{ma minda-1992-a}.
\begin{lem} \label{him-p2-lem-1.13} \cite{ma minda-1992-a}
Let $f \in \mathcal{C}(\phi)$. Then $zf''(z)/f'(z) \prec zk''(z)/k'(z)$ and $f'(z) \prec k'(z)$.
\end{lem}
\begin{lem} \label{him-p2-lem-1.14} \cite{ma minda-1992-a}
Assume $f \in \mathcal{C}(\phi)$ and $|z|=r<1$. Then 
\begin{equation} \label{him-p2-e-1.15}
-k(-r) \leq |f(z)| \leq k(r).
\end{equation}
Equality holds for some $z \neq 0$ if, and only if, f is a rotation of $k$.
\end{lem}
It is justified in \cite{ma minda-1992-a} that $-k(-r)$ is increasing in $(0,1)$ and bounded by $1$ beacuse  each $f \in \mathcal{C}(\phi)$ is normalized by $f(0)=f'(0)-1=0$. 
Therefore,  $\lim_{r \rightarrow 1} -k(-r)$ exists and denote it by $-k(-1)$.

\vspace{6mm}

Ma and Minda \cite{ma minda-1992-a} have introduced the analytic and univalent function $\phi$ with certain conditions, one of which is $\phi '(0)>0$. Recently, Kumar and Banga \cite{sivaprasadkumar-2020} have considered a non-Ma-Minda function $\Phi$, obtained by altering only one condition, namely $\Phi'(0)<0$, which is merely a rotation. Such a function $\Phi$ is named (see \cite{sivaprasadkumar-2020}) a non-Ma-Minda of type-${\bf A}$, here ${\bf A}$ signifies the condition $\Phi'(0)<0$. 
\begin{defn} \label{him-p2-def-1.2}
An analytic and univalent function $\Phi$ defined in the unit disk $\mathbb{D}$ is said to be a non-Ma-Minda of type ${\bf A}$ if it has positive real part in $\mathbb{D}$, $\Phi (\mathbb{D})$ is symmetric with respect to the real axis, starlike with respect to $\Phi(0)=1$ and $\Phi '(0)<0$. Further, it has a power series expansion of the form:
\begin{equation} \label{him-p2-e-1.16}
\Phi(z)=1+ \sum_{n=1}^{\infty} C_{n}z^{n} \quad (C_{1}<0).
\end{equation}
The class of all such functions of non-Ma-Minda of type ${\bf A}$, we shall denote it by $\mathcal{M}_{A}$.
\end{defn}  
On the similar lines of Definition \ref{him-p2-def-1.1}, Kumar and Banga has defined the classes $\mathcal{S}^{*}(\Phi)$ and $\mathcal{C}(\Phi)$. It is worth noting that by mere replacing $z$ by $-z$, {\it i.e.} $\Phi(z)=\phi(-z)$, each function in $\mathcal{M}_{A}$ reduces to its Ma-Minda counter part. For instance, $\sqrt{1+z}$ and $1-log(1-z)$ belong to the class $\mathcal{M}$ whereas  $\sqrt{1-z}$ and $1-log(1+z)$ belong to $\mathcal{M}_{A}$.

\vspace{5mm}
Observe that all the coefficients of $\phi \in \mathcal{M}$ are not necessary to be positive. For instance, some of the coefficients of $\phi (z)=z+\sqrt{1+z^{2}}= 1+z+^{2}/2-z^{4}/8+\cdot\cdot \cdot$ and $\phi(z)=1+z-z^{3}/3$ are negative. The class $\mathcal{S}^{*}(z+\sqrt{1+z^{2}})$ has been studied by Raina and Sok\'{o}\l ~ \cite{raina-2015}. Wani and Swaminathan \cite{Wani-2019} have studied the class $\mathcal{S}^{*}_{Ne}:=\mathcal{S}^{*}\left(1+z-z^{3}/3\right)$.\\

\vspace{5mm}
In 1990, Silverman and Silvia \cite{silverman-1990} introduced the following clases  $\mathcal{G}^{*}_{\alpha}$ and $\mathcal{G}_{\alpha}$. 
\begin{defn} \label{him-p2-def-1.3}
Let $\mathcal{G}^{*}_{\alpha}$ denote the class of functions $G$ analytic in $\mathbb{D}$ that satisfy the following conditions:
\begin{enumerate}
\item[(i)] $G$ is normalized by $G(0)=1$ and $G(1)=\lim_{r \rightarrow -1}=0$,

\item[(ii)] $G(\mathbb{D})$ lies in a sector with aperture $2(1-\alpha) \pi$ and vertex at the origin, and 

\item[(iii)] $G$ maps $\mathbb{D}$ univalently onto a domain that is starlike with respect to $G(1)$.
In addition, let the constant function $1$ belong to $\mathcal{G}^{*}_{\alpha}$.
\end{enumerate}
\end{defn}

\begin{defn} \label{him-p2-def-1.4}
For $\alpha$, let $\mathcal{G}_{\alpha}$ denote the class of functions, $G$ with $G(z)=1+\sum_{n=1}^{\infty} d_{n} z^{n}$, analytic and non-vanishing in $\mathbb{D}$ which  satisfy 
\begin{equation} \label{him-p2-e-1.17}
\real \left(\dfrac{zG'(z)}{G(z)} + \dfrac{(1-\alpha)(1+z)}{1-z}\right)>0
\end{equation}
for $z \in \mathbb{D}$.
\end{defn}
For $\alpha=1/2$, Robertson \cite{robertson-1981} has studied extensively the classes $\mathcal{G}^{*}_{1/2}$ and $\mathcal{G}_{1/2}$ and has conjectured that they are equal. This conjecture was proved by Lyzzaik in $1984$ (see\cite{lyzzaik-1984}). In more general, for all $\alpha$, $0\leq \alpha <1$, Silverman and Silvia have shown that $\mathcal{G}^{*}_{\alpha} = \mathcal{G}_{\alpha}$ by proving closely related property of $\mathcal{G}_{\alpha}$ with $\mathcal{S}^{*}(\alpha)$. 
\begin{lem} \label{him-p2-lem-1.18} \cite{silverman-1990}
A function $G$ is in $\mathcal{G}_{\alpha}$ if, and only if, there exists a function $s \in \mathcal{S}^{*}(\alpha)$ such that 
$$ 
G(z)=(1-z)^{2(1-\alpha)} \dfrac{s(z)}{z}.
$$
\end{lem}
The regions of variability for the class $\mathcal{G}_{\alpha}$ has been  studied extensively by Ponnusamy {\it et al.} \cite{Ponnusamy-Vasu-Vuorinen-2009}. The folowing growth theorem for the class $\mathcal{G}_{\alpha}$ has been established by  Silverman and Silvia \cite{silverman-1990}.
\begin{lem} \label{him-p2-lem-1.20} \cite{silverman-1990} 
If $g \in \mathcal{G}_{\alpha}$ then 
\begin{equation} \label{him-p2-e-1.21}
\left(\dfrac{1-r}{1+r}\right)^{2(1-\alpha)} \leq |g(z)| \leq \left(\dfrac{1+r}{1-r}\right)^{2(1-\alpha)}
\end{equation}
for $|z|=r$. Equality holds for $g(z)=((1-z)/(1+z))^{2(1-\alpha)}$ at $z=r$ and $z=-r$.
\end{lem}

In 2018, Bhowmik and Das \cite{bhowmik-2018} proved an interesting result for subordination classes. Let $f$ and $g$ be two analytic functions in $\mathbb{D}$ such that $g \prec f$. Let 
\begin{equation} \label{him-p2-e-1.22}
g(z)=\sum_{n=0}^{\infty} c_{n}z^{n}.
\end{equation}
\begin{lem} \label{him-p2-lem-1.23} \cite{bhowmik-2018}
Let $f$ and $g$ be anlytic in $\mathbb{D}$ with Taylor expansions \eqref{him-p2-e-1.1} and \eqref{him-p2-e-1.22} respectively and $g \prec f$, then 
\begin{equation} \label{him-p2-e-1.24}
\sum_{n=0}^{\infty} |c_{n}| r^{n} \leq \sum_{n=0}^{\infty} |a_{n}| r^{n}
\end{equation}
for $z|=r \leq 1/3.$
\end{lem} 

The following coefficients bounds for the class $\mathcal{S}^{*}(\alpha)$ are required to obtain the Bohr radius for the class $\mathcal{G}_{\alpha}$.
\begin{lem} \label{him-p2-lem-1.25} \cite{graham-2003}
Let $f \in \mathcal{S}^{*}(\alpha)$ be given by \eqref{him-p2-e-1.15}. Then 
$$
|a_{n}|\leq \dfrac{1}{(n-1)!} \prod_{k=0}^{n} (k-2\alpha) \quad \mbox{for} \quad n\geq 2.
$$
The equalities in above estimates are attained for $f(z)=z/(1-z)^{2(1-\alpha)}$ for $z \in \mathbb{D}$.
\end{lem}

 In this paper, we establish the Bohr phenomenon for the classes $\mathcal{S}^{*}(\phi)$ and $\mathcal{C}(\phi)$, where all the coefficients 
 of associated Ma-Minda function $\phi$ are positive such that $\phi \in H^{2}$, the Hardy class of analytic functions in $\mathbb{D}$. 
As a consequence, we obtain several important corollaries for particular choices of $\phi$.
We  also obtain the Bohr radius for the  class  $\mathcal{G}_{\alpha}$, starlike functions with respect to the bounday point.

\section{Main results}
Using Lemmas \ref{him-p2-lem-1.10}, \ref{him-p2-lem-1.11} and \ref{him-p2-lem-1.23}, we obtain Bohr radius for the class $\mathcal{S}^{*}(\phi)$.
\begin{thm} \label{him-p2-thm-2.1}
Let $f \in \mathcal{S}^{*}(\phi)$ be given by \eqref{him-p2-e-1.5} and $\phi(z)$ be given by \eqref{him-p2-e-1.6} with all $B_{n}>0, \, n \geq 1$ such that $\phi \in H^{2}$, the Hardy class of analytic functions in $\mathbb{D}$. Then 
\begin{equation} \label{him-p2-e-2.2}
|z|+ \sum_{n=0}^{\infty} |a_{n}| |z|^{n} \leq d(f(0),\partial f(\mathbb{D}))
\end{equation}
for $|z|=r \leq \min \{r_{f},1/3\}$, where $r_{f}$ is the smallest positive root of 
\begin{equation} \label{him-p2-e-2.3}
h(r)+h(-1)=0
\end{equation}
in $(0,1)$ and $h(z)$ is defined in \eqref{him-p2-e-1.9}.
\end{thm}
\begin{rem}
Since $\Phi \in \mathcal{M}_{A}$ is obtained from $\phi \in \mathcal{M}$ by mere replacing $z$ by $-z$, the image of $\mathbb{D}$ under the functions $\Phi$ and $\phi$ are identical. 
Therefore we have $\mathcal{S}^{*}(\Phi)=\mathcal{S}^{*}(\phi)$ and hence  the Bohr radius for the class $\mathcal{S}^{*}(\Phi)$ is same as that of $\mathcal{S}^{*}(\phi)$.
\end{rem}
As a consequence of Theorem \ref{him-p2-thm-2.1}, we obatin the following corollaries for particular choices of $\phi$.
\begin{cor} \label{him-p2-cor-2.4}
For $\phi(z)=\alpha +(1-\alpha)e^{z}$, we have $\mathcal{S}^{*}_{\alpha,e}:=\mathcal{S}^{*}(\alpha +(1-\alpha)e^{z})$. Let $f \in \mathcal{S}^{*}_{\alpha,e}$ be given by \eqref{him-p2-e-1.5} with $0 \leq \alpha < 0.05284$. Then the inequality \eqref{him-p2-e-2.2} is satisfied for $|z|=r\leq r_{f}$, where  $0<r_f<1/3$. The radius $r_{f}$ is the best possible. 
\end{cor}
In particular, for  $\alpha=0$ in Corollary \ref{him-p2-cor-2.4},   we obtain  the sharp Bohr radius for the class $\mathcal{S}^{*}_{e}$.

\begin{cor} \label{him-p2-cor-2.5}
For $\phi(z)=1+4z/3+2z^{2}/3$, $\mathcal{S}^{*}(\phi)$ reduces to $\mathcal{S}^{*}_{C}$. For $f \in \mathcal{S}^{*}_{C}$ of the 
form \eqref{him-p2-e-1.5}, the inequality \eqref{him-p2-e-2.2} is satisfied for $|z|=r \leq r_{f}$, where  $0<r_f<1/3$. The radius $r_{f}$ is the best possible.
\end{cor}

\begin{cor} \label{him-p2-cor-2.6}
Let $\phi$ be the rational functon such that $\phi(z)=1+(z/k) \left((k+z)/(k-z)\right)$, where $k=\sqrt{2}+1$. Then $\mathcal{S}^{*}(\phi)$ reduces to the class $\mathbb{S}^{*}_{R}$. 
Then the ineuality \eqref{him-p2-e-2.2} is satisfied for $|z|=r \leq 1/3$ for the class $\mathbb{S}^{*}_{R}$.
\end{cor}

\begin{cor} \label{him-p2-cor-2.7}
For $\phi(z)=(1+Az)/(1+Bz)$, the class $\mathcal{S}^{*}(\phi)$ reduces to the Janowski starlike class $\mathcal{S}^{*}[A,B]$. Let $f \in \mathcal{S}^{*}[A,B]$ be given by \eqref{him-p2-e-1.5} with $-1 \leq B < (1-3^{k})/(1+3^{k})<0$ and $0 \leq A \leq 1$, where $k=B/(B-A)$. Then the ineuality \eqref{him-p2-e-2.2} is satisfied for 
$|z|=r \leq r_{f}$, where  $0<r_f<1/3$. The radius $r_{f}$ is the best possible.
\end{cor}

\begin{cor} \label{him-p2-cor-2.8}
For $\phi(z)=1+z/(1-\alpha z^{2})$, the class $\mathcal{S}^{*}(\phi)$ reduces to $\mathcal{BS}^{*}(\alpha)$. Let $f \in \mathcal{BS}^{*}(\alpha)$ be given by \eqref{him-p2-e-1.5} with $0 \leq \alpha <1$. Then the inequality \eqref{him-p2-e-2.2} is satisfied for $|z|=r \leq r_{f}$, where  $0<r_f<1/3$. The constant $r_{f}$ cannot be improved.
\end{cor}

\begin{cor} \label{him-p2-cor-2.9}
Let $\phi = (1+sz)^{2}$, then $\mathcal{S}^{*}(\phi)$ reduces to the class $\mathcal{ST}_{L}(s)$. Let $f \in \mathcal{ST}_{L}(s)$ be given by \eqref{him-p2-e-1.5} with $0.444981<s \leq 1/\sqrt{2}$. Then the inequality \eqref{him-p2-e-2.2} satisfied for $|z|=r \leq r_{f}$, where  $0<r_f<1/3$. The radius $r_{f}$ is the best posible. 
\end{cor}
By using Lemmas \ref{him-p2-lem-1.14} and \ref{him-p2-lem-1.23},  we establish  Bhor phenomenon for the class $\mathcal{C}(\phi)$.
\begin{thm} \label{him-p2-thm-2.10}
Let $f \in \mathcal{C}(\phi)$ be given by \eqref{him-p2-e-1.5} and $\phi(z)$ be given by \eqref{him-p2-e-1.6} with all $B_{n}>0, \, n \geq 1$ such that $\phi \in H^{2}$, the Hardy class of analytic functions in $\mathbb{D}$. Then 
\begin{equation} \label{him-p2-e-2.11}
|z|+ \sum_{n=0}^{\infty} |a_{n}| |z|^{n} \leq d(f(0),\partial f(\mathbb{D}))
\end{equation}
for $|z|=r \leq \min \{r_{f},1/3\}$, where $r_{f}$ is the smallest positive root of 
\begin{equation} \label{him-p2-e-2.12}
k(r)+k(-1)=0
\end{equation}
in $(0,1)$ and $k(z)$ is defined in \eqref{him-p2-e-1.9}.
\end{thm}

Let $G \in \mathcal{G}_{\alpha}$ with the power series representation 
\begin{equation} \label{him-p2-e-2.13}
G(z)=1+\sum\limits_{n=1}^{\infty} d_{n} z^{n}, \quad z \in \mathbb{D}.
\end{equation}
By using Lemmas \ref{him-p2-lem-1.18},   \ref{him-p2-lem-1.20} and \ref{him-p2-lem-1.25},  we obtain the sharp  Bhor radius for the class  $ \mathcal{G}_{\alpha}$.
\begin{thm} \label{him-p2-thm-2.14}
Let $G \in \mathcal{G}_{\alpha}$ be given by \eqref{him-p2-e-2.13}. Then 
\begin{equation}
\sum\limits_{n=1}^{\infty} |d_{n}||z|^{n} \leq d(G(0), \partial G(\mathbb{D}))
\end{equation}
for $|z|=r\leq r_{f}$, where $r_{f}=(2^{1/2(1-\alpha)}-1)/(2^{1/2(1-\alpha)}+1)$. The radius $r_{f}$ is the best possible.
\end{thm}

\section{Proof of the main results}

\begin{pf} [{\bf Proof of Theorem   \ref{him-p2-thm-2.1}}]
Let $f \in \mathcal{S}^{*}(\phi)$ then in view of Lemma \ref{him-p2-lem-1.11}, we have 
\begin{equation} \label{him-p2-e-3.1}
|f(z)|\geq -h(-r)  \quad \mbox{for } \quad |z|<1.
\end{equation}
By taking $\liminf$ as $|z|=r\rightarrow 1$ on both the sides of \eqref{him-p2-e-3.1}, we obtain 
\begin{equation} \label{him-p2-e-3.2}
\liminf \limits_{|z|\rightarrow 1} |f(z)| \geq -h(-1).
\end{equation}
The Euclidean distance between $f(0)$ and the boundary of $f(\mathbb{D})$ is given by 
\begin{equation} \label{him-p2-e-3.3}
d(f(0), \partial f(\mathbb{D}))= \liminf \limits_{|z|\rightarrow 1} |f(z)-f(0)|.
\end{equation}
Since $f(0)=0$ from \eqref{him-p2-e-3.2} and \eqref{him-p2-e-3.3}, we obtain 
\begin{equation} \label{him-p2-e-3.4}
d(f(0), \partial f(\mathbb{D})) \geq -h(-1).
\end{equation} 
It is known that $\phi \in H^{2}$ if, and only if, 
\begin{equation} \label{him-p2-e-3.5}
\sum\limits_{n=1}^{\infty} |B_{n}|^{2}< + \infty.
\end{equation}
From \eqref{him-p2-e-3.5} we have
\begin{equation} \label{him-p2-e-3.6}
\sum\limits_{n=1}^{\infty} \dfrac{B_{n}}{n} \leq \sqrt{\sum\limits_{n=1}^{\infty} |B_{n}|^{2}} \, \sqrt{\sum\limits_{n=1}^{\infty}\dfrac{1}{n^{2}}} < + \infty.
\end{equation}
 Let 
\begin{equation} \label{him-p2-e-3.7}
h(z)=z+ \sum\limits_{n=2}^{\infty} b_{n} \,z^{n} \quad \mbox{for } \quad z\in \mathbb{D}.
\end{equation}
In view of \eqref{him-p2-e-1.9}, we obtain $zh'(z)/h(z)=\phi (z)$ and a simple computation shows that 
\begin{equation} \label{him-p2-e-3.8}
h(z)= z \exp \left(\int\limits_{0}^{z} \dfrac{\phi(t)-1}{t}\, dt\right),
\end{equation}
where integgration is along the linear segment joing $0$ to $z\in \mathbb{D}$. 
From \eqref{him-p2-e-1.9} and \eqref{him-p2-e-3.8}, we btain 
\begin{equation} \label{him-p2-e-3.9}
h(z)=z\, \exp \, \left(\sum\limits_{n=1}^{\infty} B_{n} \dfrac{z^{n}}{n}\right).
\end{equation}
Let $H:[0,1] \rightarrow \mathbb{R}$ be defined by 
$$
H(r)=h(r)+h(-1).
$$
From \eqref{him-p2-e-3.6}, \eqref{him-p2-e-3.8} and \eqref{him-p2-e-3.9}, it is easy to see that $H$ is continuous in $[0,1]$ and differentiable in $(0,1)$. Note that $H(0)=h(-1)<0$ and $H(1)=h(1)+h(-1)$. Clearly, 
$$
h(1)=\exp\, \left(\sum\limits_{n=1}^{\infty} \dfrac{B_{n}}{n}\right) \quad \mbox{and} \quad h(-1)= \exp\, \left(\sum\limits_{n=1}^{\infty}  \dfrac{(-1)^{n} B_{n}}{n}\right).
$$
Since $B_{n}>0$ for $n\geq 1$, 
$$
\sum\limits_{n=1}^{\infty} \dfrac{B_{n}}{n}> \sum\limits_{n=1}^{\infty}  \dfrac{(-1)^{n} B_{n}}{n}.
$$
Therefore, 
$$
h(1)=\exp\, \left(\sum\limits_{n=1}^{\infty} \dfrac{B_{n}}{n}\right)>\exp\, \left(\sum\limits_{n=1}^{\infty}  \dfrac{(-1)^{n} B_{n}}{n}\right)=-h(-1)
$$
which implies that $H(1)>0$. Since $H(0)<0$ and $H(1)>0$, by the intermediate value property, $H$ has one real root in $(0,1)$. Let $r_{f}$ be the smallest positive root of $H$ in (0,1). Thus $H(r_{f})=0$, which is equivalent to 
\begin{equation} \label{him-p2-e-3.10}
h(r_{f})=-h(-1).
\end{equation}
 From \eqref{him-p2-e-3.7} and \eqref{him-p2-e-3.9}, we obtain
\begin{equation} \label{him-p2-e-3.11}
h(z)=z+ \sum\limits_{n=2}^{\infty} b_{n} \,z^{n} = z\, \exp \, \left(\sum\limits_{n=1}^{\infty} B_{n} \dfrac{z^{n}}{n}\right)  \quad \mbox{for } \quad z\in \mathbb{D}.
\end{equation}
Therefore,
\begin{align} \label{him-p2-e-3.12}
|z|+\sum\limits_{n=2}^{\infty} |b_{n}||z|^{n} 
& = |z|\exp \, \left(\sum\limits_{n=1}^{\infty} |B_{n}| \dfrac{|z|^{n}}{n}\right) \\[5mm] \nonumber
 & = r \exp\, \left(\sum\limits_{n=1}^{\infty} B_{n} \dfrac{r^{n}}{n}\right) 
 \\ \nonumber
  &= h(r). \end{align}
  In view of Lemma \ref{him-p2-lem-1.23} and \eqref{him-p2-e-3.11}, for $|z|=r$, we obtain 
\begin{equation} \label{him-p2-e-3.13}
r+ \sum\limits_{n=2}^{\infty} |a_{n}|r^{n} \leq r+\sum\limits_{n=2}^{\infty} |b_{n}|r^{n} = h(r)  \quad \mbox{for} \quad r \leq \dfrac{1}{3}.
\end{equation}
For $r \leq r_{f}$, we have 
\begin{equation} \label{him-p2-e-3.14}
h(r)\leq h(r_{f})=-h(-1).
\end{equation}
Combining \eqref{him-p2-e-3.4}, \eqref{him-p2-e-3.13} and \eqref{him-p2-e-3.14}, we conclude that 
$$
r+\sum\limits_{n=2}^{\infty} |a_{n}|r^{n} \leq -h(-1)  \leq d(f(0),\partial f(\mathbb{D}))
$$ 
for $r \leq \min \{r_{f},1/3\}$.  
If $r_{f}$ lies in $(0,1/3]$, then $r_{f}$ is the best possible. Let $0<r_{f}\leq 1/3$. To show the sharpness of $r_{f}$, we choose $f$ to be a suitable rotation $h_{\epsilon}$ of the function $h$ {\it i.e.} $f=h_{\epsilon}$ which is obviously belongs to $\mathcal{S}^{*}(\phi)$. Since equality in Lemma \ref{him-p2-lem-1.11} occurs for a suitable rotation of $h$, we have 
\begin{equation} \label{him-p2-e-3.15}
d(h_{\epsilon}(0),\partial h_{\epsilon}(\mathbb{D}))=-h_{\epsilon}(-1).
\end{equation}
A simple computation using \eqref{him-p2-e-3.10} and \eqref{him-p2-e-3.15}, for $f=h_{\epsilon}$ and  $|z|=r_{f}$ shows that
$$
|z|+ \sum_{n=0}^{\infty} |a_{n}| |z|^{n}=h_{\epsilon}(r)= -h_{\epsilon}(-1)=d(h_{\epsilon}(0),\partial h_{\epsilon}(\mathbb{D})).
$$
Therefore, the radius $r_{f} \in (0,1/3)$ is the best possible. This completes the proof.
\end{pf}

\begin{pf}[\bf Proof of Corollary \ref{him-p2-cor-2.4}]
Given that $\phi(z)=\alpha +(1-\alpha) e^{z}$ and each coefficient of $\phi(z)$ is strictly positive for $0 \leq \alpha <1$. By the similar lines of argument as in the proof of Theorem \ref{him-p2-thm-2.1}, consider $H(r)=h(r)+h(-1)$. Note that 
$$
h\left(\dfrac{1}{3}\right)=\dfrac{1}{3} \exp\, \left((1-\alpha) \int\limits_{0}^{\frac{1}{3}}\left(\dfrac{-1+e^{t}}{t}\right)\, dt\right) \approx \dfrac{1}{3} (1.43807)^{1-\alpha}
$$
and 
$$
h(-1)=-\exp\, \left((1-\alpha) \int\limits_{0}^{-1}\left(\dfrac{-1+e^{t}}{t}\right)\, dt\right)\approx -(0.450859463)^{1-\alpha}.
$$
Therefore, using Mathematica we can see that $H(1/3)=h(1/3)+h(-1)>0$ if, and only if, $0 \leq \alpha < 0.05284$. On the other hand, $H(0)=h(-1)<0$. 
Therefore we conclude that $H$ has a  root in $(0,1)$. Since $r_{f}$ is the smallest root of $H$ in $(0,1)$, we have $r_{f}<1/3$. Hence in view of Theorem \ref{him-p2-thm-2.1}, $r_{f}$ is the best possible.
\end{pf}

\begin{pf}[\bf Proof of Corollary \ref{him-p2-cor-2.5}]
Let $\phi(z)=1+4z/3+2z^{2}/3$. A simple computation shows that 
$$
h(r)=r \exp\, \left(\dfrac{4}{3}r+\dfrac{r^{2}}{3}\right).
$$ 
Note that $h(1/3) \approx 0.539490$ and $h(-1)\approx -0.367879441$. Let
$$
 H(r)=h(r)+h(-1).
 $$
  Then $H(0) \approx -0.3678799441<0$ and $H(1/3)>0$. Therefore, $H$ has a root in $(0,1/3)$. Let $r_{f}$ be the smallest root in $(0,1/3)$. In view of Theorem  \ref{him-p2-thm-2.1}, $r_{f}$ is the best possible.
\end{pf}

\begin{pf}[\bf Proof of Corollary \ref{him-p2-cor-2.6}]
For $\phi (z)=1+(z/k) \left((k+z)/(k-z)\right)$, where $k=\sqrt{2}+1$, a simple computation using \eqref{him-p2-e-3.8} shows that 
$$
h(r)=\dfrac{r}{e^{r}} \left(\dfrac{k}{k-r}\right)^{2k}.
$$	
Note that $h(-1)\approx-0.5099807$ and $h(1/3)\approx 0.489391446$. Hence $H(1/3)=h(1/3)+h(-1)<0$ and $H(0)<0$. Therefore $H$ has no root in $(0,1/3)$ and $r_{f}>1/3$. Hence the inequality $\eqref{him-p2-e-2.2}$ is satisfied for $r \leq 1/3$. 
\end{pf}

\begin{pf}[\bf Proof of Corollary \ref{him-p2-cor-2.7}]
Let $\phi (z)=(1+Az)/(1+Bz)$ with $-1 \leq B<A\leq 1$. It is easy to see that the coefficients of $\phi(z)$ are all positive when $-1 \leq B \leq 0$. A simple computation using \eqref{him-p2-e-3.8} yields 
$$
h(r)=r \left(1+Br \right)^{\left(\frac{A-B}{B}\right)}.
$$ 
Let $H(r)= h(r)+h(-1)$. Then $H(0)<0$ and 
$$
H\left(\dfrac{1}{3}\right)= h\left(\dfrac{1}{3}\right) + h(-1)=\dfrac{1}{3}\left(1+\dfrac{B}{3}\right)^{\left(\frac{A-B}{B}\right)}-\left(1-B\right)^{\left(\frac{A-B}{B}\right)}.
$$
A simple computation shows that $H(1/3)>0$ when $-1 \leq B \leq (1-3^{k})/(1+3^{k})$, where $k=(B-A)/B$. Therefore, $H$ has a root in $(0,1)$. Let the smallest root of $H$ in $(0,1/3)$ be $r_{f}$. In view of Theorem \ref{him-p2-thm-2.1}, we conclude that for $r \leq r_{f} $, the inequality $\eqref{him-p2-e-2.2}$ is satisfied and the radius $r_{f}$ is the best possible.  
\end{pf}

\begin{pf} [\bf Proof of Corollary \ref{him-p2-cor-2.8}]
Here $\phi(z)=1+z/(1-\alpha z^{2})$. Using \eqref{him-p2-e-3.8}, we obtain 
$$
h(r)=r \left(\dfrac{1+\sqrt{\alpha}\, r}{1-\sqrt{\alpha}\, r}\right)^{\frac{1}{2\sqrt{\alpha}}}.
$$
Clearly, $h(-1)=-\left(\dfrac{1+\sqrt{\alpha}}{1+\sqrt{\alpha}}\right)^{\frac{1}{2\sqrt{\alpha}}}$ and 
$h\left(\dfrac{1}{3}\right)=\left(\dfrac{1+\dfrac{\sqrt{\alpha}}{3}}{1+\dfrac{\sqrt{\alpha}}{3}}\right)^{\frac{1}{2\sqrt{\alpha}}}.$	\\
Let $H(r)= h(r)+h(-1)$. Then $H(0)<0$ and $H(1/3)=h(1/3)+h(-1)>0$ for $0 \leq \alpha <1$. Therefore, $H$ has a root in $(0,1/3)$ and choose $r_{f}$ to be 
the smallest root in $(0,1/3)$. By Theorem \ref{him-p2-thm-2.1}, the radius $r_{f}$ is the best possible.
\end{pf}

\begin{pf}[\bf Proof of Corollary \ref{him-p2-cor-2.9}]
Let $\phi (z)= (1+sz)^{2}$ with $0<s \leq 1/\sqrt{2}$. Using \eqref{him-p2-e-3.8}, we obtain 
$$
h(r)=r\exp \left(s\left(2r+\dfrac{sr^{2}}{2}\right)\right).
$$
Let $H(r)=h(r)+h(-1)$. Note that $H(0)<0$ and 
$$
H\left(\dfrac{1}{3}\right)=\dfrac{1}{3} \exp \left(s\left(\dfrac{s+12}{18}\right)\right)- \exp \left(s\left(-2+\dfrac{s}{2}\right)\right)>0,
$$
if $0.444981< s \leq 1/\sqrt{2}$. Therefore, $H$ has a real root in $(0,1/3)$. Let $r_{f}$ be the smallest root in $(0,1/3)$. In view of Theorem \ref{him-p2-thm-2.1}, the radius $r_{f}$ is the best possible.
\end{pf}

\begin{pf}[\bf Proof of the Theorem \ref{him-p2-thm-2.10}]
Let $f \in \mathcal{C}(\phi)$ then from Lemma \ref{him-p2-lem-1.14}, it is evident that  
\begin{equation} \label{him-p2-e-3.15-a}
d(f(0), \partial f(\mathbb{D}))= \liminf \limits_{|z|\rightarrow 1} |f(z)-f(0)|\geq -k(-1).
\end{equation}
It is known that $f \in \mathcal{C}(\phi)$ if, and only if, $zf' \in \mathcal{S}^{*}(\phi)$. Therefore, there exists $g \in  \mathcal{S}^{*}(\phi)$ such that $zf'(z)=g(z)$ for $z \in \mathbb{D}$. This relation gives 
\begin{equation} \label{him-p2-e-3.16}
f(z)=\int\limits_{0}^{z} \dfrac{g(t)}{t} \,\, dt.
\end{equation} 
Let $f(z)=z+ \sum_{n=2}^{\infty}a_{n}z^{n}$ and $g(z)=z+ \sum_{n=2}^{\infty} b_{n}z^{n}$. Then from \eqref{him-p2-e-3.16}, we obtain
\begin{equation} \label{him-p2-e-3.17}
z+ \sum\limits_{n=2}^{\infty}a_{n}z^{n}= z+ \sum\limits_{n=2}^{\infty} \dfrac{b_{n}}{n} z^{n}.
\end{equation}
Note that the majorant series for $g$, we have $M_{g}(r):=r+\sum_{n=2}^{\infty} |b_{n}|r^{n}$. Then 
\begin{align} \label{him-p2-e-3.18}
\int\limits_{0}^{r} \dfrac{M_{g}(t)}{t} \,\, dt 
&=r+ \sum\limits_{n=2}^{\infty} |b_{n}| \dfrac{r^{n}}{n} 
\\ \nonumber & = r+ \sum\limits_{n=2}^{\infty} |a_{n}|r^{n}
\\ \nonumber & = M_{f}(r).
\end{align}
Since we have $zk'(z)=h(z)$, then going by same lines of argument as in \eqref{him-p2-e-3.18}, we obtain 
\begin{equation} \label{him-p2-e-3.19}
M_{k}(r)=\int\limits_{0}^{r} \dfrac{M_{g}(t)}{t} \, dt.
\end{equation} 
In view of Lemma \ref{him-p2-lem-1.23}, we have 
\begin{equation} \label{him-p2-e-3.20}
M_{g}(r) \leq M_{h}(r) \quad \mbox{for} \quad r\leq \dfrac{1}{3}.
\end{equation}
Using \eqref{him-p2-e-3.18}, \eqref{him-p2-e-3.19} and \eqref{him-p2-e-3.20}, for $r \leq 1/3$, we obtain 
\begin{equation} \label{him-p2-e-3.21}
M_{f}(r) \leq \int\limits_{0}^{r} \dfrac{M_{h}(t)}{t} \, dt= M_{k}(r)=k(r).
\end{equation}
By the same lines of argument as in Theorem \ref{him-p2-thm-2.1}, we can show that $H_{1}: [0,1] \rightarrow \mathbb{R}$ defined by $H_{1}(r)=k(r)+k(-1)$ has a real root in $(0,1)$. Let $r_{f}$ be the smallest root of $H_{1}$ in $(0,1)$ and $H_{1}(r_{f})=0$, which follows that $k(r_{f})=-k(-1)$. Therefore, for $r \leq \min \{r_{f},1/3\} $, we obtain 
$$
M_{f}(r) \leq k(r)  \leq -k(-1)  \leq d(f(0), \partial f(\mathbb{D})).
$$
Let $0<r_{f}\leq 1/3$. To show the sharpness of $r_{f}$, we choose $f$ to be a suitable rotation $k_{\epsilon}$ of the function $k$ {\it i.e.} $f=k_{\epsilon}$ which is obviously belongs to $\mathcal{C}(\phi)$. Since the equality in Lemma \ref{him-p2-lem-1.14} occurs for a suitable rotation of $k$, we obtain
\begin{equation} \label{him-p2-e-3.22}
d(k_{\epsilon}(0),\partial k_{\epsilon}(\mathbb{D}))=-k_{\epsilon}(-1).
\end{equation}
A simple computation using \eqref{him-p2-e-3.15-a}, for $f=k_{\epsilon}$ and  $|z|=r_{f}$, shows that 
$$
|z|+ \sum_{n=0}^{\infty} |a_{n}| |z|^{n}=k_{\epsilon}(r)= -k_{\epsilon}(-1)=d(k_{\epsilon}(0),\partial k_{\epsilon}(\mathbb{D})).
$$
Therefore, the radius $r_{f} \in (0,1/3)$ is the best possible. This completes the proof.
\end{pf}

\begin{pf}[\bf Proof of the Theorem \ref{him-p2-thm-2.14}]
Let $G \in \mathcal{G}_{\alpha}$ be given by \eqref{him-p2-e-2.13}. Then from Lemma \ref{him-p2-lem-1.20}, it is evident that 
\begin{equation} \label{him-p2-e-3.23}
d(G(0), \partial G(\mathbb{D}))= \liminf \limits_{|z|\rightarrow 1} |G(z)-G(0)| \geq 1.
\end{equation}
In view of Lemma \ref{him-p2-lem-1.18}, we have the following relation
\begin{equation} \label{him-p2-e-3.24}
G(z)=(1-z)^{2(1-\alpha)} \dfrac{s(z)}{z}.
\end{equation}
For $s\in \mathcal{S}^{*}(\alpha)$, using Lemma \ref{him-p2-lem-1.25} we obtain 
\begin{equation} \label{him-p2-e-3.25}
M_{s}(r) \leq \dfrac{r}{(1-r)^{2(1-\alpha)}}
\end{equation}
for $0<r<1$, where $M_{s}(r)$ is the associated majorant series of $s$.
It is known that for $f$ and $g$ be two analytic functions in $\mathbb{D}$, 
\begin{equation} \label{him-p2-e-3.26}
M_{fg} \leq M_{f} M_{g}
\end{equation}
for $0<r<1$, where $M_{fg}$, $M_{f}$ and $M_{g}$ are the associated majorant series of $fg$, $f$ and $g$ respectively. 
In view of \eqref{him-p2-e-3.24}, \eqref{him-p2-e-3.25} and \eqref{him-p2-e-3.26}, for $|z|=r$, we obtain 
\begin{align} \label{him-p2-e-3.27}
M_{G}(r) &= 1+ \sum\limits_{n=1}^{\infty} |d_{n}||z|^{n}
\\ \nonumber & \leq (1+|z|)^{2(1-\alpha)} \dfrac{M_{s}(r)}{r}
\\ \nonumber & \leq \dfrac{(1+r)^{2(1-\alpha)}}{(1-r)^{2(1-\alpha)}}.
\end{align}
Therefore from \eqref{him-p2-e-3.27}, we obtain 
\begin{equation} \label{him-p2-e-3.27-a}
\sum\limits_{n=1}^{\infty} |d_{n}||z|^{n} \leq \dfrac{(1+r)^{2(1-\alpha)}}{(1-r)^{2(1-\alpha)}} -1.
\end{equation}
The right hand side of \eqref{him-p2-e-3.27-a} is less than or equal to $1 \leq d(G(0), \partial G(\mathbb{D}))$ if 
$$
\dfrac{(1+r)^{2(1-\alpha)}}{(1-r)^{2(1-\alpha)}}\leq 2
$$	
{\it i.e.} for $r \leq r_{G}:=(2^{1/2(1-\alpha)}-1)/(2^{1/2(1-\alpha)}+1)$.
To show that $r_{G}$ is the best possible, we consider the function $G_{\alpha}:\mathbb{D}\rightarrow \mathbb{C}$ defined by 
$$
G_{\alpha}(z)=\left(\dfrac{1-z}{1+z}\right)^{2(1-\alpha)}.
$$
From Lemma \ref{him-p2-lem-1.20}, for $G=G_{\alpha}$, we have 
$$
d(G(0), \partial G(\mathbb{D}))=1.
$$
For $G=G_{\alpha}$ and $|z|=r_{G}$, a simple computation shows that 
$$
\sum\limits_{n=1}^{\infty} |d_{n}||z|^{n} =\left(\dfrac{1+r_{G}}{1-r_{G}}\right)^{2(1-\alpha)}-1=1=d(G(0), \partial G(\mathbb{D})).
$$
This shows that the radius $r_{f}$ is the best posible. This completes the proof.
\end{pf}

\noindent\textbf{Acknowledgement:}  The first author thank SERB-MATRICS and the second author thank CSIR for their support.


\begin{thebibliography}{99}
	
	
	
\bibitem{Abu-2010} {\sc Y. Abu-Muhanna},  Bohr's phenomenon in subordination and bounded harmonic classes, {\it Complex Var. Elliptic Equ.} {\bf  55} (2010), 1071--1078.
	
		
	\bibitem{abu-2011} {\sc Y. Abu-Muhanna} and  {\sc R. M. Ali}, Bohr's phenomenon for analytic functions into the exterior of a compact convex body, {\it J. Math. Anal. Appl.} {\bf  379}  (2011), 512--517.
		
	
	\bibitem{abu-2013} {\sc Y.  Abu Muhanna} and {\sc R. M.  Ali}, Bohr's phenomenon for analytic functions and the hyperbolic metric, {\it Math. Nachr.}  {\bf 286}  (2013), 1059--1065.
	
		
	\bibitem{abu-2014} {\sc Y.Abu Muhanna, R. M. Ali, Z. C. Ng} and  {\sc S. F. M Hasni}, Bohr radius for subordinating families of analytic functions and bounded harmonic mappings, 
	{\it J. Math. Anal.Appl.} {\bf 420} (2014), 124--136.
	
		
	\bibitem{aizn-2000} {\sc L.Aizenberg}, Multidimensional analogues of Bohr's theorem on power series, \textit{Proc. Amer. Math. Soc.} {\bf 128} (2000), 1147--1155.
	
		
	\bibitem{aizenberg-2001} {\sc L. Aizenberg, A. Aytuna}  and {\sc P. Djakov}, Generalization of theorem on Bohr for bases in spaces of holomorphic functions of several complex variables, 
	{\it J. Math. Anal.Appl.} {\bf  258} (2001), 429--447.
	
	
	\bibitem{aizn-2007} {\sc L.Aizenberg}, Generalization of results about the Bohr radius for power series, {\it Stud. Math.}  {\bf 180}  (2007), 161--168.  
	
		
	%\bibitem{Ali 1995} {\sc R. M .Ali}, {\sc S. Ponnusamy} and {\sc V. Singh}, Starlikeness of functions satisfying a differential inequality,  {\it Ann. Polon. Math.} {\bf  61} (1995), 135--140.
		
	\bibitem{Ali-2017} {\sc R. M. Ali, R.W. Barnard} and {\sc  A.Yu. Solynin}, A note on Bohr's phenomenon for power series, {\it J. Math. Anal.Appl.} {\bf 449} (2017), 154-167.
	
	\bibitem{Ali-2019} {\sc R. M. Ali, N. K.Jain} and {\sc V. Ravichandran}, Bohr radius for classes of analytic functions, {\it Results Math.} {\bf 74} (2019).
		
	\bibitem{alkhaleefah-2019} {\sc S. A. Alkhaleefah, I.R. Kayumov} and {\sc S. Ponnusamy}, On the Bohr inequality with a fixed zero coefficient, {\it Proc. Amer. Math. Soc.} {\bf 147} (2019), 5263--5274.
	
\bibitem{Himadri-Vasu-P1} {\sc Vasudevarao Allu} and {\sc Himadri Halder},  Bhor phenomenon for certain subclasses of Harmonic Mappings, arXiv:2006.11622, 2020
	
	\bibitem{bene-2004} {\sc C. B{\rm $\acute{E}$}n$\acute{E}$teau, A. Dahlner} and {\sc D. Khavinson}, Remarks on the Bohr phenomenon, {\it  Comput. Methods Funct. Theory}  {\bf 4}  (2004), 1--19.
		
	\bibitem{bhowmik-2018} {\sc B. Bhowmik} and {\sc N. Das},  Bohr phenomenon for subordinating families of certain univalent functions,  {\it J. Math. Anal. Appl.}  {\bf 462} (2018), 1087--1098.
	
	\bibitem{boas-1997} {\sc H.P. Boas} and {\sc D. Khavinson}, Bohr's power series theorem in several variables, {\it Proc. Amer. Math. Soc}  {\bf 125} (1997), 2975--2979.
	
	\bibitem{boas-2000} {\sc H.P. Boas}, Majorant Series,  {\it J. Korean Math. Soc.}  {\bf 37} (2000), 321--337.
			
	\bibitem{Bohr-1914} {\sc H. Bohr}, A theorem concerning power series,  {\it Proc. Lond. Math. Soc}. s2-13 (1914), 1--5.

	\bibitem{evdoridis-2019} {\sc S. Evdoridis} and {\sc S. Ponnusamy}, Improved Bohr's inequality for locally univalent harmonic mappings, {\it  Indag. Math. (N.S.) } {\bf 30} (2019), 201--213.
	
	\bibitem{duren-1983} {\sc P. L. Duren}, {\it Univalent Functions} (Grundlehren der mathematischen Wisseenschaften 259, New York, Berlin, Heidelberg, Tokyo) Springer-Verlag, 1983.
	
	\bibitem{goodman-1991-a} {\sc A. W. Goodman}, On uniformly convex functions, {\it Annl. Polon. Math.}, {\bf 56} (1991), 87--92.
	
	\bibitem{goodman-1991-b} {\sc A. W. Goodman}, On uniformly starlike functions, {\it J. Math. Anal. Appl.}, {\bf 155} (1991), 364--370.
	
	\bibitem{graham-2003} {\sc I. Graham} and {\sc G. Kohr}, Geometric Function Theory in One and Higher Dimensions, Monographs and Textbooks in Pure and Applied Mathematics, {\bf 255}, Marcel Dekker, Inc., New York, 2003. 
	\bibitem{kargar-2019} {\sc R. Kargar}, {\sc A. Ebadian} and {\sc J. Sokol}, On Both lemniscate and starlike functions, {\it Anal. Math. Phys.} {\bf 9} (2019), 143--154. 
	
%	\bibitem{kanas-2019} {\sc S. Kanas}, {\sc V. S. Masih} and {\sc A. Ebadian}, Relations of a planar domains bounded by hyperbolas with families of holomorphic functions, {\it J. Inequal. Appl.} {\bf 246} (2019), 1--14.
	\bibitem{kayumov-2018-a} {\sc I.R. Kayumov, S. Ponnusamy} and {\sc N. Shakirov}, Bohr radius for localy unovalent harmonic mappings, {\it  Math. Nachr.} {\bf  291} (2018), 1757--1768.
	
	\bibitem{kayumov-2018-b} {\sc I.R. Kayumov} and {\sc S. Ponnusamy}, Bohr's inequalities for the analytic functions with lacunary series and harmonic functions, 
	{\it J. Math. Anal. Appl.}  {\bf 465} (2018), 857--871.
	
	\bibitem{kayumov-2019} {\sc I.R. Kayumov} and {\sc S. Ponnusamy}, On a powered Bohr inequality, {\it  Ann. Acad. Sci. Fenn. Math.} {\bf 44} (2019), 301--310.
	
	\bibitem{kayumova-2018} {\sc A. Kayumova, I. R. Kayumov} and {\sc S. Ponnusamy}, Bohr's inequality for harmonic mappings and beyond,
	 {\it Mathematics and Computing}, 245--256, Commun. Comput. Inf. Sci., 834, Springer, Singapore, 2018.
	
	\bibitem{khatter-2019} {\sc K.Khatter}, {\sc V. Ravichandran} and {\sc S. Sivaprasad Kumar}, Starlike functions associated with exponential function and the lemniscate of Bernoulli, 
	{\it Rev. R. Acad. Cienc. Exactas F$\acute{i}$s. Nat. Ser. A Mat. RACSAM} {\bf 113} (2019), 233--253.
	
	\bibitem{kumar-2016} {\sc S. Kumar} and {\sc V. Ravichandran}, A subclass of starlike functions associated with a rational function, {\it Southeast Asian Bull. Math.} {\bf 40} (2016), 199--212.
	
	\bibitem{sivaprasadkumar-2020} {\sc S. Sivaprasad Kumar} and {\sc S. Banga}, On a special type of Non-Ma-Minda function, arXiv:2006.02111v1, 2020.
	
	\bibitem{Liu-2019} {\sc Z. Liu} and {\sc S. Ponnusamy}, Bohr radius for subordination and $k$-quasiconformal harmonic mappings, {\it  Bull. Malays. Math. Sci. Soc.} {\bf 42} (2019), 2151--2168.
	
	\bibitem{lyzzaik-1984} {\sc A. Lyzzail}, On a conjecture of M. S. Robertson, {\it Proc. Amer. Math. Soc.} {\bf 91} (1984), 108--110. 
	
	\bibitem{ma minda-1992-a} {\sc  W.C. Ma} and {\sc  D. Minda}, A unified treatment of some special classes of univalent functions, in {\it Proceedings of the Conference on Complex Analysis}(Tianjin, 1992), 157--169, Conf. Proc. Lecture Notes Anal., I, Int. Press, Cambridge.
	
	\bibitem{ma minda-1992-b} {\sc  W.C. Ma} and {\sc  D. Minda}, Uniformly convex functions, {\it Ann. Polon. Math.} {\bf 57} (1992), 165--175.
	
	\bibitem{ma minda-1993} {\sc  W.C. Ma} and {\sc  D. Minda}, Uniformly convex functions II, {\it Ann. Polon. Math.} {\bf 58} (1993), 275--285.
	
	\bibitem{masih-2020} {\sc V. S. Masih} and {\sc S. Kanas}, Subclasses of Starlike and Convex functions associated with associated with the Limaçon Domain, {\it Symmetry} {\bf 12} (2020).
	\bibitem{mendiratta-2015} {\sc R. Mendiratta}, {\sc S. Nagpal} and {\sc V. Ravichandran}, On a subclass of strongly starlike functions associated with exponential function, {\it Bull. Malays. Math. Sci. Soc.} {\bf 38} (2015), 365--386.
	\bibitem{paulsen-2002} {\sc Vern I. Paulsen, Gelu Popescu} and {\sc Dinesh Singh}, On Bohr's inequality, {\it Proc. Lond. Math. Soc.} s3-85 (2002), 493--512.
	
	\bibitem{Ponnusamy-Vasu-Vuorinen-2009}	{\sc S. Ponnusamy,  A. Vasudevarao} and {\sc  M. Vuorinen}, Region of variability for spirallike functions with respect to a boundary point, 
         {\it Colloq. Math.} {\bf 116} (2009), 31--46.


	\bibitem{raina-2015} {\sc R. K.Raina} and {\sc J. Sokol}, Some properties related to a certain class of starlike functions, {\it C. R. Math. Acad. Sci. Paris} {\bf 353} (2015), 973--978.
	
	\bibitem{robertson-1981} {M. S. Robertson}, Univalent functions with respect to a boundary point, {\it J. Math. Anal. Appl.} {\bf 81} (1981), 327--345.
	\bibitem{ronning-1991} {\sc F. Ronning}, On starlike functions associated with parabolic regions, {\it Anna. Univ. Mariae Curie-Sklodowska Sect. A}, {\bf 45} (1991), 117--122.
	
	\bibitem{ronning-1993} {\sc F. Ronning}, Uniformly convex functions and a corresponding class of starlike functions, {\it Proc. Amer. Math. Soc.} {\bf 118} (1993), 189--196.
	
	\bibitem{sharma-2016} {\sc K. Sharma}, {\sc N. K. Jain} and {\sc V. Ravichandran}, Starlike functions associated with cardiod, {\it Afr. Mat.} {\bf 27} (2016), 923--939.
	
	\bibitem{sidon-1927} {\sc S. Sidon}, Uber einen satz von Hernn Bohr, {\it Math. Zeit.}  {\bf 26} (1927),  731-732.
	
	\bibitem{silverman-1990} {\sc H. Silverman} and {\sc E.M. Silvia}, Subclasses of univalent functions starlike with respect to a boundary point, {\it Houston J. Math.} {\bf 16} (1990), 289--299. 
	
	\bibitem{vasu-book} {\sc Derek K. Thomas}, {\sc Nikola Tuneski} and {\sc Allu Vasudevarao}, Univalent functions. A primer, De Gruyter Studies in Mathematics, {\bf 69}. De Gruyter, Berlin, 2018.
	
	\bibitem{tomic-1962} {\sc M. Tomic}, Sur un theoreme de H. Bohr,  {\it Math. Scand.}  {\bf 11} (1962), 103--106.
	
	\bibitem{Wani-2019} {\sc L. A. Wani} and {\sc A. Swaminathan}, Starlike and convex functions associated with nephroid domain, {\it Bull. Malays. Math. Sci. Soc.} doi:10.1007/s40840-020-00935-6. 2020
\end{thebibliography}
\end{document}